# A Geometric Proof that *e* is Irrational and a New Measure of its Irrationality

## Jonathan Sondow

**1. INTRODUCTION.** While there exist geometric proofs of irrationality for √2 [**2**], [**27**], no such proof for $e$, $\pi$, or $\ln 2$ seems to be known. In section 2 we use a geometric construction to prove that $e$ is irrational. (For other proofs, see [**1**, pp. 27-28], [**3**, p. 352], [**6**], [**10**, pp. 78-79], [**15**, p. 301], [**16**], [**17**, p. 11], [**19**], [**20**], and [**21**, p. 302].) The proof leads in section 3 to a new measure of irrationality for $e$, that is, a lower bound on the distance from $e$ to a given rational number, as a function of its denominator. A connection with the greatest prime factor of a number is discussed in section 4. In section 5 we compare the new irrationality measure for $e$ with a known one, and state a number-theoretic conjecture that implies the known measure is almost always stronger. The new measure is applied in section 6 to prove a special case of a result from [**24**], leading to another conjecture. Finally, in section 7 we recall a theorem of G. Cantor that can be proved by a similar construction.

**2. PROOF.** The irrationality of $e$ is a consequence of the following construction of a nested sequence of closed intervals $I_n$. Let $I_1 = [2,3]$. Proceeding inductively, divide the interval $I_{n-1}$ into $n$ ($\geq 2$) equal subintervals, and let the second one be $I_n$ (see Figure 1). For example, $I_2 = \left[\frac{5}{2!}, \frac{6}{2!}\right]$, $I_3 = \left[\frac{16}{3!}, \frac{17}{3!}\right]$, and $I_4 = \left[\frac{65}{4!}, \frac{66}{4!}\right]$.

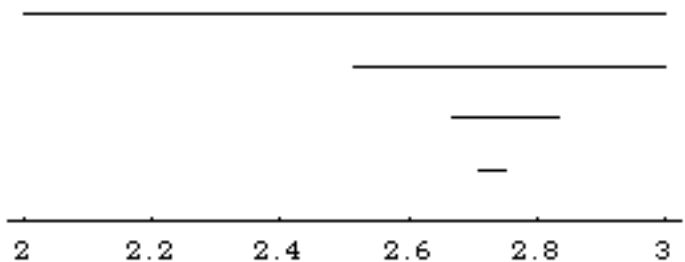

**Figure 1.** The intervals $I_1, I_2, I_3, I_4$.

The intersection

$$\bigcap_{n=1}^{\infty} I_n = \{e\} \qquad (1)$$

is then the geometric equivalent of the summation (see the Addendum)

$$\sum_{n=0}^{\infty} \frac{1}{n!} = e. \qquad (2)$$

When $n > 1$ the interval $I_{n+1}$ lies strictly between the endpoints of $I_n$, which are $\frac{a}{n!}$ and $\frac{a+1}{n!}$ for some integer $a = a(n)$. It follows that the point of intersection (1) is not a fraction with denominator $n!$ for any $n \geq 1$. Since a rational number $p/q$ with $q > 0$ can be written

$$\frac{p}{q} = \frac{p \cdot (q-1)!}{q!}, \tag{3}$$

we conclude that $e$ is irrational. •

**Question.** The nested intervals $I_n$ intersect in a number—let's call it $b$. It is seen by the Taylor series (2) for $e$ that $b = e$. Using only standard facts about the natural logarithm (including its definition as an integral), but *not* using any series representation for log, can one see directly from the given construction that $\log b = 1$?

**3. A NEW IRRATIONALITY MEASURE FOR $e$.** As a bonus, the proof leads to the following measure of irrationality for $e$.

**Theorem 1.** *For all integers p and q with $q > 1$*

$$\left| e - \frac{p}{q} \right| > \frac{1}{(S(q)+1)!}, \tag{4}$$

*where $S(q)$ is the smallest positive integer such that $S(q)!$ is a multiple of $q$.*

For instance, $S(q) = q$ if $1 \leq q \leq 5$, while $S(6) = 3$. In 1918 A. J. Kempner [13] used the prime factorization of $q$ to give the first algorithm for computing

$$S(q) = \min\{k > 0 : q \mid k!\} \tag{5}$$

(the so-called Smarandache function [28]). We do not use the algorithm in this note.

*Proof of Theorem 1.* For $n > 1$ the left endpoint of $I_n$ is the closest fraction to $e$ with denominator not exceeding $n!$. Since $e$ lies in the interior of the second subinterval of $I_n$,

$$\left| e - \frac{m}{n!} \right| > \frac{1}{(n+1)!} \tag{6}$$

for any integer $m$. Now given integers $p$ and $q$ with $q > 1$, let $m = p \cdot S(q)!/q$ and $n = S(q)$. In view of (5), $m$ and $n$ are integers. Moreover,

$$\frac{p}{q} = \frac{p \cdot S(q)!/q}{S(q)!} = \frac{m}{n!}. \tag{7}$$

Therefore, (6) implies (4). •

As an example, take $q$ to be a prime. Clearly, $S(q) = q$. In this case, (4) is the (very weak) inequality

$$\left| e - \frac{p}{q} \right| > \frac{1}{(q+1)!}. \tag{8}$$

In fact, (4) implies that (8) holds for *any* integer $q$ larger than 1, because $S(q) \leq q$ always holds. But (4) is an improvement of (8), just as (7) is a refinement of (3).

Theorem 1 would be false if we replaced the denominator on the right side of (4) with a smaller factorial. To see this, let $p/q$ be an endpoint of $I_n$, which has length $\frac{1}{n!}$. If we take $q = n!$, then since evidently

$$S(n!) = n \tag{9}$$

and $e$ lies in the interior of $I_n$,

$$\left| e - \frac{p}{q} \right| < \frac{1}{S(q)!}. \tag{10}$$

(If $q < n!$, then (10) still holds, since $n > 2$, so $p/q$ is not an endpoint of $I_{n-1}$, hence $S(q) = n$.)

**4. THE LARGEST PRIME FACTOR OF $q$.** For $q \geq 2$ let $P(q)$ denote the largest prime factor of $q$. Note that $S(q) \geq P(q)$. Also, $S(q) = P(q)$ if and only if $S(q)$ is prime. (If $S(q)$ were prime but greater than $P(q)$, then since $q$ divides $S(q)!$, it would also divide $(S(q) - 1)!$, contradicting the minimality of $S(q)$.)

P. Erdős and I. Kastanas [**9**] observed that

$$S(q) = P(q) \quad \text{(almost all } q\text{)}. \tag{11}$$

(Recall that a claim $C_q$ is true *for almost all $q$* if the counting function $N(x) = \#\{q \leq x : C_q \text{ is false}\}$ satisfies the asymptotic condition $N(x)/x \to 0$ as $x \to \infty$.) It follows that Theorem 1 implies an irrationality measure for $e$ involving the simpler function $P(q)$.

**Corollary 1.** *For almost all $q$, the following inequality holds with any integer $p$:*

$$\left| e - \frac{p}{q} \right| > \frac{1}{(P(q)+1)!}. \tag{12}$$

When $q$ is a factorial, the statement is more definite.

**Corollary 2.** *Fix $q = n! > 1$. Then (12) holds for all $p$ if and only if $n$ is prime.*

*Proof.* If $n$ is prime, then $P(q) = n$, so (4) and (9) imply (12) for all $p$. Conversely, if $n$ is composite, then $P(q) < n$, and (10) shows that (12) fails for certain $p$.  •

Thus when $q > 1$ is a factorial, (12) is true for all $p$ if and only if $S(q) = P(q)$. To illustrate this, take $\frac{p}{q} = \frac{65}{4!}$ to be the left endpoint of $I_4$. Then $P(q) = 3 < 4 = S(q)$, and (12) does not hold, although of course (4) does:

$$0.00833\ldots = \frac{1}{5!} < \left|e - \frac{65}{24}\right| = 0.00994\ldots < \frac{1}{4!} = 0.04166\ldots.$$

**5. A KNOWN IRRATIONALITY MEASURE FOR $e$.** The following measure of irrationality for $e$ is well known: *given any $\varepsilon > 0$ there exists a positive constant $q(\varepsilon)$ such that*

$$\left|e - \frac{p}{q}\right| > \frac{1}{q^{2+\varepsilon}} \tag{13}$$

*for all $p$ and $q$ with $q \geq q(\varepsilon)$.* This follows easily from the continued fraction expansion of $e$. (See, for example, [**23**]. For sharper inequalities than (13), see [**3**, Corollary 11.1], [**4**], [**7**], [**10**, pp. 112-113], and especially the elegant [**26**].)

Presumably, (13) is usually stronger than (4). We state this more precisely, and in a number-theoretic way that does not involve $e$.

**Conjecture 1.** *The inequality $q^2 < S(q)!$ holds for almost all $q$. Equivalently, $q^2 < P(q)!$ for almost all $q$.*

(The equivalence follows from (11).) This is no doubt true; the only thing lacking is a proof. (Compare [**12**], where A. Ivić proves an asymptotic formula for the counting function $N(x) = \#\{q \leq x : P(q) < S(q)\}$ and surveys earlier work, including [**9**].)

Conjecture 1 implies that (13) is almost always a better measure of irrationality for $e$ than those in Theorem 1 and Corollary 1. On the other hand, Theorem 1 applies to all $q > 1$. Moreover, (4) is stronger than (13) for certain $q$. For example, let $q = n!$ once more. Then (4) and (9) give (6), which is stronger than (13) if $n > 2$, since

$$(n+1)! < (n!)^2 \quad (n \geq 3). \tag{14}$$

**6. PARTIAL SUMS VS. CONVERGENTS.** Theorem 1 yields other results on rational approximations to $e$ [**24**]. One is that *for almost all $n$, the $n$-th partial sum $s_n$ of series* (2) *for $e$ is* not *a convergent to the simple continued fraction for $e$.* Here $s_0 = 1$ and $s_n$ is the left endpoint of $I_n$ for $n \geq 1$. (In 1840 J. Liouville [**14**] used the partial sums of the Taylor series for $e^2$ and $e^{-2}$ to prove that the equation $ae^2 + be^{-2} = c$ is impossible if $a$, $b$, and $c$ are integers with $a \neq 0$. In particular, $e^4$ is irrational.)

Let $q_n$ be the denominator of $s_n$ in lowest terms. When $q_n = n!$ (see [**22**, sequence A102470]), the result is more definite, and the proof is easy.

**Corollary 3.** *If $q_n = n!$ with $n \geq 3$, then $s_n$ cannot be a convergent to $e$.*

*Proof.* Use (4), (9), (14), and the fact that every convergent satisfies the reverse of inequality (13) with $\varepsilon = 0$ [**10**, p. 24], [**17**, p. 61]. •

When $q_n < n!$ (for example, $q_{19} = 19!/4000$—see [**22**, sequence A093101]), another argument is required, and we can only prove the assertion for almost all $n$. However, numerical evidence suggests that much more is true.

**Conjecture 2.** *Only two partial sums of series* (2) *for e are convergents to e, namely, $s_1 = 2$ and $s_3 = 8/3$.*

**7. CANTOR'S THEOREM.** A generalization of the construction in section 2 can be used to prove the following result of Cantor [**5**].

**Theorem 2.** *Let $a_0, a_1, \ldots$ and $b_1, b_2, \ldots$ be integers satisfying the inequalities $b_n \geq 2$ and $0 \leq a_n \leq b_n - 1$ for all $n \geq 1$. Assume that each prime divides infinitely many of the $b_n$. Then the sum of the convergent series*

$$a_0 + \frac{a_1}{b_1} + \frac{a_2}{b_1 b_2} + \frac{a_3}{b_1 b_2 b_3} + \cdots$$

*is irrational if and only if both $a_n > 0$ and $a_n < b_n - 1$ hold infinitely often.*

For example, series (2) for $e$ and all subseries (such as $\Sigma_{n \geq 0} \frac{1}{(2n)!} = \cosh 1$ and $\Sigma_{n \geq 0} \frac{1}{(2n+1)!} = \sinh 1$) are irrational, but the sum $\Sigma_{n \geq 1} \frac{n-1}{n!} = 1$ is rational.

An exposition of the "if" part of Cantor's theorem is given in [**17**, pp. 7-11]. For extensions of the theorem, see [**8**], [**11**], [**18**], and [**25**].

**ADDENDUM.** Here are some details on why the nested closed intervals $I_n$ constructed in section 2 have intersection $e$. Recall that $I_1 = [2, 3]$, and that for $n \geq 2$ we get $I_n$ from $I_{n-1}$ by cutting it into $n$ equal subintervals and taking the second one. The left-hand endpoints of $I_1, I_2, I_3, \ldots$ are $2, 2 + \frac{1}{2!}, 2 + \frac{1}{2!} + \frac{1}{3!}, \ldots$, which are also partial sums of the series (2) for $e$. Since the endpoints approach the intersection of the intervals, whose lengths tend to zero, the intersection is the single point $e$.

**ACKNOWLEDGMENTS.** Stefan Krämer pointed out the lack of geometric proofs of irrationality. The referee suggested a version of the question in section 2. Yann Bugeaud and Wadim Zudilin supplied references on the known irrationality measures for $e$. Aleksandar Ivić commented on Conjecture 1. Kyle Schalm did calculations [**24**] on Conjecture 2, and Yuri Nesterenko related it to Liouville's proof. I am grateful to them all.

*209 West 97th Street, New York, NY 10025*
*jsondow@alumni.princeton.edu*